\documentclass{amsart}

\newtheorem{Theorem}{Theorem}

\usepackage{graphics}
\begin{document}
\title{On the cordial deficiency of complete multipartite graphs}
\author{Adrian Riskin}
\address{Department of mathematics\\
Mary Baldwin College\\
Staunton, Virginia 24401\\
USA\\}
\email{ariskin@mbc.edu}
\thanks{I would like to thank Sam Patteson and his colleagues
 for providing an ideal working environment
for the writing of this paper.}
\subjclass{05C78}
\keywords{Cordial graph}

\begin{abstract}
We calculate the cordial edge deficiencies of the complete multipartite graphs and find an 
upper bound for their cordial vertex deficiencies.  We also give conditions under which the 
tensor product of two cordial graphs is cordial.
\end{abstract}

\maketitle

\baselineskip=20pt

\section{Introduction and definitions}

Cahit [2] introduced cordial graph labelings as a generalization of both graceful and 
harmonious labelings.  See Gallian [3] for a comprehensive bibliography on the topic.
We allow graphs to have multiple edges but not loops.  A \textit{binary labeling} of a 
graph $G$ is a map $f:V(G) \rightarrow \{0,1\}$.  Two real numbers are said to be
\textit{roughly equal} if $0 \leq |x-y| \leq 1$. A binary labeling of $G$ is \textit{friendly}
provided that $|f^{-1}(0)|$ is roughly equal to $|f^{-1}(1)|$.  A binary labeling $f$ of $G$
induces a labeling $f_{e}:E(G) \rightarrow \{0,1\}$ by $f_{e}(uv)=f(u)+f(v)$ where
$uv \in E(G)$ and the sum is calculated modulo 2. A friendly labeling $f$ of $G$ is called 
\textit{cordial} when $|f_{e}^{-1}(0)|$ is roughly equal to $|f_{e}^{-1}(1)|$, and a graph 
$G$ is called \textit{cordial} if if has a cordial labeling.

In [6] I introduced a new measure of the degree to which a noncordial graph fails to be cordial, 
inspired by Kotzig's and Rosa's notion of \textit{edge-magic deficiency} [5], and studied another
such measure, first defined in [1]. Note first 
that every friendly labeling $f$ of a graph can be made into a cordial labeling of an augmented
graph $G^{\prime}$  by adding no more than $|f_{e}^{-1}(0) - f_{e}^{-1}(1)|-1$ edges between 
appropriate pairs of vertices so that $|f_{e}^{-1}(0)|$ becomes roughly equal to $|f_{e}^{-1}(1)|$.
The minimum number of edges, taken over all friendly labelings of $G$, which it is necessary to add
in order that $G^{\prime}$ become cordial is called the \textit{cordial edge deficiency} of $G$, 
denoted by $\mathrm{ced}(G)$.  This is essentially the same concept as 
the \textit{index of cordiality} introduced in [1].  If it is possible to find a binary labeling $f$ of $G$ so that
$|f_{e}^{-1}(0)|$ and $|f_{e}^{-1}(1)|$ are roughly equal, then it is possible to make $f$ into 
a cordial labeling of an augmented graph $G^{\prime}$ by adding no more than 
$|f^{-1}(0) - f^{-1}(1)|-1$ isolated vertices labeled in such a way as to make $f$ into a friendly
labeling.  The minimum number of vertices, taken over all such binary labelings of $G$, which it
is necessary to add in order to make $G^{\prime}$ cordial is called the \textit{cordial vertex
deficiency} of $G$, denoted by $\mathrm{cvd(G)}$.  If there are no such binary labelings of $G$
we call $G$ \textit{strictly noncordial} and write $\mathrm{cvd}(G) = \infty$.  If $G$ is a graph
and $f$ is a binary labeling of $G$ then the \textit{cordial deficit} of the pair $(G,f)$ is
$||f_{e}^{-1}(0)| - |f_{e}^{-1}(1)||$.

Finally, we will be studying the cordiality of the tensor products of graphs.  If $G$ and $H$ are
graphs, then the tensor product $G \times H$ is the graph whose vertex set  
is the cartesian product of the vertex sets of $G$ and $H$, namely $V(G) \times V(H)$, in which
$\left( (u_{1}, u_{2}), (v_{1}, v_{2})\right) \in E(G \times H)$ if and only if $u_{1}v_{1} \in E(G))$ 
and $u_{2}v_{2} \in E(H)$.  The tensor product is also known as the weak product and as the 
categorical product of $G$ and $H$.  Considerable effort has focused on tensor products due to 
Hedetniemi's conjecture, for a useful survey of which see  [7]¥.

\section{Comments on Lee's and Liu's constructions}

In 1991 Lee and Liu [4] published the following theorem:

\begin{Theorem}
Let H be a graph with an even number of edges and a cordial labeling such that the vertices of H can be divided
into $\ell$ parts $H_{1}, H_{{2}}, \dots , H_{\ell}$ each consisting of an equal number of vertices labeled 0 and vertices
labeled 1.  Let $G$ be any graph and $G_{1}, G_{2}, \dots, G_{\ell}$ be any $\ell$ subsets of the vertices of G.
Let (G,H) be the graph which is the disjoint union of $G$ and $H$ augmented by edges joining every vertex in $G_{i}$ to 
every vertex in $H_{1}$ for $1 \leq \ell$.  Then $G$ is cordial if and only if $(G,H)$ is.
\end{Theorem}

They provide an explicit proof that the cordiality of $G$ implies the cordiality of $(G,H)$ and state that ``the converse 
can be proved in the same way.''  Unfortunately this is not the case for the theorem as it is stated.   In fact the 
converse as stated is false.  The theorem would in fact be true if it were required that the restriction of the cordial
labeling of $(G,H)$ to the vertices of $H$ is a cordial labeling of $H$  This need not be true for every cordial labeling
of $(G,H)$ as we will show below with a counterexample.  It is extremely important to note that the only problem with 
this theorem is its statement. In fact, Lee and Liu use the correct version throughout the paper, so that not only are the 
other theorems in the paper correct, but the proofs are in fact correct as well.  In order for the statement of the 
theorem to be correct, the last sentence must be replaced with:
\begin{quote}
Then $G$ is cordial if and only if $(G,H)$ has a cordial labeling which, when restricted to $H$ yields an equal number of
vertices labeled 0 and vertices labeled 1 in each $H_{i}$
\end{quote}

\noindent \textbf{Counterexample 1:}
Let $H=C_{4}$.  Then (i) $H$ is cordial,  (ii) $H$ has an even number of edges, and (iii) $H$ has a cordial labeling such
that the vertices of $H$ can be divided into two parts $H_{1}$ and $H_{2}$ such that the number of vertices labeled 0 is
equal to the number of vertices labeled 1 (see Figure 1).  

\begin {figure}[h]
\begin{center}
\caption{}
\includegraphics{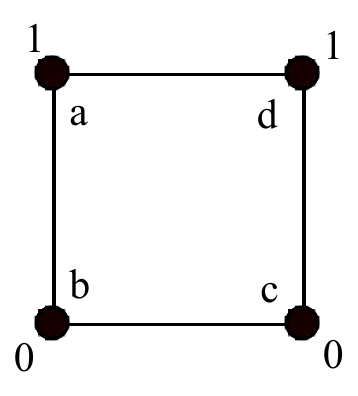}
\end{center}
\end{figure}

\noindent Note that $H_{1}=\{a,b\}$ and $H_{2}=\{c,d\}$.  Now let $G=K_{4}$.
Cahit [2] showed that $K_{n}$ is cordial if and only if $n \leq 3$.  We label the vertices of $K_{4}$ as shown in 
Figure 2.

\begin {figure}[h]
\begin{center}
\caption{}
\includegraphics{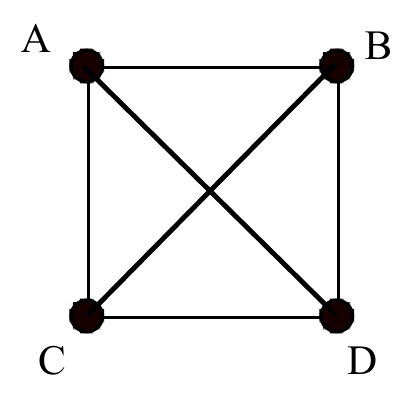}
\end{center}
\end{figure}

\noindent Let $G_{1}=\{A,B\}$ and $G_{2}=\{B,C\}$.  It is easily checked that the following is a cordial labeling of
$(G,H)$:

\begin{center}
$\begin{matrix}
a & b & c & d & A & B & C & D\\
0 & 1 & 1 & 1 & 0 & 0 & 1 & 0
\end{matrix}$
\end{center}

\section{Cordial edge-deficiency of complete multipartite graphs}

One of Lee and Liu's most interesting results is:

\begin{Theorem}
A complete $k$-partite graph is cordial if and only if the number of parts with an odd number of
vertices is at most three.
\end{Theorem}

\noindent We generalize this using the concept of cordial edge-deficiency thus:

\begin{Theorem}
Let $G$ be a complete multipartite graph with $k$ odd parts.  Then ced($G$)$=max
\left \{0,\left \lfloor \frac
{k}{2}\right \rfloor-1  \right \}$¥
\end{Theorem}

\noindent \textbf{Proof:} Let $f:V(G) \rightarrow \{0,1\}$ be a friendly labeling of $G$.
Let $E_{1}, \dots , E_{j}$ be the even parts of $G$. Suppose that an even part $E_{i}$ is 
assigned fewer zeros than ones by the labeling $f$.  Because $f$ is a friendly labeling
there must be another part of $G$ which has more zeros than ones.  There are two cases to
consider, depending on whether this other part is even or odd.  In the first case, suppose that
$E_{\ell}$ has more zeros than ones. Let $E_{i}$ have $z_{i}$ zeros and $E_{\ell}$ have
$z_{\ell}$ zeros.  Furthermore, assume that $|E_{i}|=2m_{i}$ and that $|E_{\ell}|=2m_{\ell}$.  
Note that $E_{i}$ has $2m_{i}-z_{i}$ ones and that $E_{\ell}$ has $2m_{\ell}-z_{\ell}$ ones.
If we switch a zero label from $E_{i}$ with a zero label from $E_{\ell}$ we produce a new 
friendly labeling $g:V(G)\rightarrow\{0,1\}$ in which there is no change in the number of edges 
labeled zero nor in the number of edges labeled one except possibly amongst the edges with one
end in $E_{i}$ and the other end in $E_{\ell}$.  Note that

\begin{equation*}
|f_{e*}^{-1}(0)|=z_{i}z_{\ell}+(2m_{i}-z_{i})(2m_{\ell}-z_{\ell})
\end{equation*}

\noindent and that

\begin{equation*}
|f_{e*}^{-1}(1)|=z_{i}(2m_{\ell}-z_{\ell}) + z_{\ell}(2m_{i}-z_{i})
\end{equation*}

\noindent where $f_{e*}$ represents $f_{e}$ restricted to the complete bipartite subgraph
of G generated by $E_{i}$ and $E_{\ell}$.  The difference between these is therefore 

\begin{equation*}
||f_{e*}^{-1}(0)|-|f_{e*}^{-1}(1)|| = 4(m_{i}-z_{i})(z_{\ell}-m_{\ell})
\end{equation*}

After the switch is completed, $E_{i}$ has $z_{i}+1$ zeros and $(2m_{i}-z_{i}-1)$ ones and 
$E_{\ell}$ has $z_{\ell} -1$ zeros and $2m_{\ell}-z_{\ell}+1$ ones.  Hence, as above, 

\begin{equation*}
|g_{e*}^{-1}(0)| = (z_{i}+1)(z_{\ell}-1)+(2m_{i}-z_{i}-1)(2m_{\ell}-z_{\ell}+1)
\end{equation*}

and

\begin{equation*}
|g_{e*}^{-1}(1)| = (z_{i}+1)(2m_{\ell}-z_{\ell}+1)+(z_{\ell}-1)(2m_{i}-z_{i}+1)
\end{equation*}

so that

\begin{equation*}
||g_{e*}^{-1}(0)|-|g_{e*}^{-1}(1)|| = 4|m_{i}-(z_{i}+1)||m_{\ell}-(z_{\ell}-1)|
\end{equation*}

Furthermore

\begin{equation*}
0 \leq m_{i}-(z_{i}+1) < m_{i}-z_{i}
\end{equation*}

and

\begin{equation*}
0 \leq (z_{\ell}-1)-m_{\ell} < z_{\ell}-m_{\ell}
\end{equation*}

and therefore

\begin{equation*}
||g_{e*}^{-1}(0)|-|g_{e*}^{-1}(1)|| < ||f_{e*}^{-1}(0)|-|f_{e*}^{-1}(1)||
\end{equation*}

Similar calculations show the analogous results for the other cases where one or both parts are odd.
It follows by induction that the cordial deficit for a friendly labeling is minimized when the 
numbers of zeros and of ones in each part are roughly equal.

Now, we may assume without loss of generality that if $k$ is odd, $|f^{-1}(0)| = |f^{-1}(1)|-1$.
Note also that the cordial deficit of the labeling is equal to the sum of the cordial deficits
of pairs of parts over all such pairs.  Furthermore, if one or both of the parts in a pair has
an even number of vertices then the cordial deficit of that pair is zero.  Hence we need only
calculate the sum over all pairs of odd parts.  If $k$ is odd, then since 
$|f^{-1}(0)| = |f^{-1}(1)|-1$ there are $\frac{k+1}{2}$ odd parts which have one more one than 
zero, and $\frac{k-1}{2}$ parts which have one more zero than one.  This makes the net cordial
deficit

\begin{equation*}
\left | \begin{pmatrix} \frac{k-1}{2}\\
2 \end{pmatrix} + \begin{pmatrix} \frac{k+1}{2}\\
2 \end{pmatrix} - \frac{k-1}{2} \frac{k+1}{2} \right |
\end{equation*}

and if $k$ is even the net cordial deficit is

\begin{equation*}
\left | 2 \begin{pmatrix} \frac{k}{2}\\
2 \end{pmatrix} - \left (\frac{k}{2}\right)^{2} \right|
\end{equation*}

In either case, ced($G$) = $\left \lfloor \frac{k}{2} \right \rfloor -1$ \hfill $\square$

The calculation of the cordial vertex deficiency of the complete multipartite graphs seems to 
be a more difficult problem.  I was able, however, to obtain an upper bound which applies in 
certain cases.  The following theorem from [6] is necessary for the proof:

\begin{Theorem}
The cordial vertex deficiency of $K_{n}$ is $j-1$ if $n=j^{2}+\delta$, where $\delta \in \{-2, 0, 2\}$.
Otherwise $K_{n}$ is strictly noncordial.
\end{Theorem}

\begin{Theorem}
If $G$ is a complete multipartite graph with $n \geq 1$ odd parts and $n=j^{2}+\delta$ where
$\delta \in \{-2,0,2\}$ then cvd$(G) \leq j-1$.
\end{Theorem}

\noindent \textbf{Proof:} Let $v_{i}$ be a single vertex from the $i^{th}$ odd part of $G$ for
$1 \leq i \leq n$.  Label the vertices in each of the even parts with half zeros and half ones.
Label the vertices in each of the odd parts, omitting $v_{i}$, with half zeros and half ones.
We now have a labeling of all but $n$ of the vertices of $G$ which has the same number of
vertices labeled zero as are labeled one and with all the edges in the subgraph $H$ isomorphic
to $K_{n}$ induced by the $v_{i}$'s labeled half with zeros and half with ones.  Now apply 
the previous theorem to $H$. \hfill $\square$

\section{Cordiality of tensor products}

Note that if $G_{1}$ is connected, simple, and bipartite, and $G_{2}$ has $q$ edges, then the tensor
product $G_{1} \times G_{2}$ is decomposable into $2q$ edge-disjoint isomorphs of $G_{1}$.  Also, if
$G_{1}$ and $G_{2}$ are cordially labeled by friendly labelings $f$ and $g$ respectively, then the
\textit{induced labeling} of the tensor product is obtained by labeling $(u,v) \in V(G_{1} \times G_{2})$
by $f(u)+g(v)$, where the sum is calculated modulo 2.

\begin{Theorem}
Let $G_{1}$ and $G_{2}$ be cordially labeled simple graphs such that $G_{1}$ is connected, bipartite, and 
has an even number of edges.  Then $G_{1} \times G_{2}$ is cordially labeled by the induced vertex
labeling.
\end{Theorem}

\noindent \textbf{Proof:} Each of the $2q$ isomorphs of $G_{1}$ is generated from a directed edge of $G_{2}$ by
using it to generate a labeling of the unique bipartition of $G_{1}$.  The induced labeling restricted to a 
particular isomorph of $G_{1}$ is obtained from the cordial labeling of $G_{1}$ by adding the label of the edge
which is generating the isomorph to the labeling of each edge of the isomorph as determined by the isomorphism
between it and cordially labeled $G_{1}$.  Since $G_{1}$ has an even number of edges, half of the edges of the 
isomorph end up labeled zero and the other half end up labeled one.  Since the isomorphs of $G_{1}$ partition 
the edges of $G_{1} \times G_{2}$, the tensor product itself ends up cordially labeled. \hfill $\square$

\bigskip

\noindent It is fairly easy to find cordial labelings of such tensor products when $G_{1}$ has an odd number of
edges.  However, it is also easy to show that in none of these cases is the induced labeling cordial.  However, I 
feel it is worth a conjecture to the effect that the requirement that $G_{1}$ have an even number of edges can
be dropped from the statement of Theorem 6.

\section{References}
\begin{enumerate}
\item Boxwala, S. and Limaye, N.B.  On the cordiality of elongated $t$-plys.
JCMCC 52(2005) 181-221.

\item Cahit, I.; Cordial graphs: a weaker version of graceful and harmonious
graphs.  Ars Combin. 23(1987) 201-207.

\item Gallian, J. A.; A dynamic survey of graph labeling.  Electronic J.
Combin. DS6. http://www.combinatorics.org/Surveys/index.html

\item Lee, S.M. and Liu, A. A construction of cordial graphs from smaller cordial graphs.
Ars Combin. 32 (1991), 209--214.

\item Kotzig, A. and Rosa, A.; Magic valuations of finite graphs.  Canad.
Math. Bull. 13(1970) 451-461.

\item Riskin, A. Cordial Deficiency.  To appear in the Bulletin of the Malaysian Journal of
Mathematics.  \{arXiv:math.CO/0610760\}

\item Zhu, X. A survey on Hedetniemi's conjecture. Taiwanese J. Math. 2(1998) 1-24. 
\end{enumerate}

\end{document}